%Format: AMSppt
%%% the version mailed to Ian on January 22

\input amstex
\documentstyle{amsppt}

\TagsOnRight
\NoBlackBoxes
\magnification=1200
\pageheight{9 true in}
\pagewidth{6.05 truein}

\font \smallcaps=cmcsc8 at 8pt
\font \sm=cmr8 at 8pt
\font \mayor=cmbx12 at 12pt
\font\smit=cmti8 at 8pt
\font\smb=cmbx8 at 8pt

\newbox\bigstrutbox
\setbox\bigstrutbox=\hbox{\vrule height9.5pt depth3.5pt width0pt}
\def\bigstrut{\relax\ifmmode\copy\bigstrutbox\else\unhcopy\bigstrutbox\fi}

\def\K{K(n)^*}

\def\ftwo{{\bold F}_2}
\def\fthree{{\bold F}_3}
\def\ffive{{\bold F}_5}
\def\fp{{\bold F}_p}
\def\Im{\hbox{\rm Im}}
\def\zz{{\bold Z}}
\def\cc{{\bold C}}
\def\qq{{\bold Q}}
\def\Hom{\hbox{\rm Hom}}
\def\Aut{\hbox{\rm Aut}}
\def\End{\hbox{\rm End}}
\def\Soc{\hbox{\rm Soc}}
\def\fs{{+_F}}
\def\Gal{\hbox{\rm Gal}}

\def\alp{{\bf 1}}
\def\bru{{\bf 2}}
\def\cur{{\bf 3}}
\def\hkr{{\bf 4}}
\def\jrh{{\bf 5}}
\def\kri{{\bf 6}}
\def\kuh{{\bf 7}}
\def\kui{{\bf 8}}
\def\rawi{{\bf 9}}
\def\gap{{\bf 10}}
\def\bjoe{{\bf 11}}

\def\teyai{{\bf 13}}
\def\teyaii{{\bf 14}}
\def\wil{{\bf 15}}
\def\wilk{{\bf 16}}
\def\yag{{\bf 17}}

\rightheadtext{On the $GL(V)$--module structure of $\K(BV)$}
\leftheadtext{Ian J. Leary and Bj\"orn Schuster}

\centerline{\mayor On the $GL(V)$--module structure of $K(n)^*(BV)$}
\bigskip
\centerline{\smc By IAN J. LEARY}
\medskip
\centerline{\it Max-Planck Institut f\"ur Mathematik, 53225 Bonn, Germany}
\medskip
\centerline{\smc and BJ\"ORN SCHUSTER}
\medskip
\centerline{\it Centre de Recerca Matem\`atica, E--08193 Bellaterra (Barcelona),
Spain}
\bigskip
\centerline{({\it January 22, 1996})}
\vskip 1 cm
\document
\head {\it Abstract}\endhead
We study the question of whether the Morava K-theory of the 
classifying space of an elementary abelian group $V$ is a 
permutation module (in either of two distinct senses, defined below) for the 
automorphism group of $V$.  We use Brauer characters and computer 
calculations.  Our algorithm for finding permutation submodules of
modules for $p$-groups may be of independent interest.  

\head 1. {\it Introduction} \endhead

Let $p$ be a prime, 
let $\K$ denote the $n$th Morava K-theory, $V$ an elementary abelian
$p$-group, or equivalently an $\fp$-vector space, and
$GL(V)$ the group of automorphisms of $V$. Then $GL(V)$ acts naturally on
the classifying space $BV$ of $V$ and hence on $h^*(BV)$ for any cohomology
theory $h$.  In the case when $h=K(n)$, $\K(BV)$ is a finitely
generated free module over the coefficient ring $\K$ whose structure 
is known [\rawi], and it is natural to ask what may be said about its
structure as a module for the group ring $\K[GL(V)]$.  The Morava
K-theory of arbitrary finite groups is not known, and there is no
direct construction of Morava K-theory itself.  We hope that a better
understanding of $\K(BV)$ may lead to progress with these questions.

For any ring $R$, and finite group $G$, we say that an 
$R$-free $R[G]$-module $M$ is a permutation module if there is an
$R$-basis for $M$ which is permuted by the action of $G$.  Call 
such an $R$-basis a permutation basis for $M$.  
If $S$ is a $G$-set, write $R[S]$ for the permutation module with
permutation basis $S$.  If $M$ is a
graded module for the graded ring $R[G]$ (where elements of $G$ are
given grading zero), we call $M$ a graded permutation module if 
it is a permutation module with a permutation basis consisting of 
homogeneous elements.  

If $M$ is a graded module for $K(n)^*[G]$, and $G_p$ is a Sylow 
$p$-subgroup of $G$, the following four conditions on $M$ are 
progressively weaker, in the sense that each is implied by the 
previous one.  

(1) $M$ is a graded permutation module; 

(2) $M$ is a permutation module; 

(3) $M$ is a direct summand of a permutation module; 

(4) as a $K(n)^*[G_p]$-module, $M$ is a graded permutation module.  

\noindent
The implications $(1)\implies (2)\implies (3)$ are obvious,
and hold for $M$ a graded $R[G]$-module for any $R$.  The implication 
$(3)\implies (4)$ is explained below in Sections 2~and~7.  

One might hope for $K(n)^*(BV)$ to satisfy condition (1) for
$G=GL(V)$, i.e., for $K(n)^*(BV)$ to be a graded permutation module
for the group ring $\K[GL(V)]$.  
This would be useful for the following reason:
The ordinary cohomology of a group with coefficients in 
any permutation module is determined by the Eckmann-Shapiro lemma. 
Hence if $\K(BV)$ is a graded permutation module for $GL(V)$, 
and $H$ is a group expressed as an extension with kernel 
$V$, the $E_2$ page of the Atiyah-Hirzebruch spectral sequence 
converging to $\K(BH)$ is easily computable.  Even condition (4) would
be very useful, as it would facilitate the computation of the
$E_2$-page of the Atiyah-Hirzebruch spectral sequence for any
$p$-group $H$ expressed as an extension with kernel $V$.  

The recent work of I.~Kriz on Morava K-theory, including his dramatic
discovery of a 3-group $G$ such that $K(2)^*(BG)$ is not concentrated
in even degrees, has emphasised the importance of studying the
$\Aut(H)$-module structure of $\K(BH)$ [\kri].  For example, Kriz has
shown that for any prime $p$ and any cyclic $p$-subgroup $C$ of
$GL(V)$, $\K(BV)$ is a (graded) permutation module for $C$.  He uses 
this result to deduce that for $p$ odd and $G$ a split extension with
kernel $V$ and quotient $C$, $\K(BG)$ is concentrated in even degrees.
N.~Yagita has another proof of this result [\yag].  

It should be noted that if the dimension of $V$ is at least three,
there are infinitely many indecomposable graded $\K[GL(V)]$-modules,
of which only finitely many occur as summands of modules satisfying 
condition (4), which suggests that a \lq random' module
will not satisfy any of the conditions.  On the other hand, work of 
Hopkins, Kuhn and Ravenel [\hkr] shows that for certain generalized
cohomology theories $h^*$, $h^*(BV)$ is a permutation module for 
$h^*[GL(V)]$.  Amongst these $h^*$ are theories closely related to 
$\K$, albeit that their coefficient rings are torsion-free and contain
an inverse for $p$.  

A result due to Kuhn [\kuh] shows that $\K(BV)$ has the same Brauer
character as the permutation module $\K[\Hom(V,(\fp)^n)]$ for $GL(V)$.  
We shall show however that in general $\K(BV)$ does not have the same 
Brauer character as a graded permutation module for $\K[GL(V)]$.  
Note that Brauer characters give no information whatsoever concerning 
the structure of $\K(BV)$ as a module for the Sylow $p$-subgroup of
$GL(V)$.  We give an algorithm to determine, for any $p$-group $G$, 
whether an $\fp[G]$-module is a permutation module, and use this 
algorithm and computer calculations to determine in some cases 
whether $\K(BV)$ satisfies condition (4) above.  

The main results of this paper are summarised in the following 
four statements.  Before making them, we fix some notation.  

\proclaim{\smc Definition} Throughout the paper, let $p$ be a prime, 
$\K$ the $n$th Morava K-theory (at the prime $p$), and let $V$ be a 
vector space over the field of $p$ elements of dimension $d$.  
Let $GL(V)$ act on the right of $V$, which will have the advantage
that the modules we consider will be left modules.  
Let $U(V)$ be a Sylow $p$-subgroup of $GL(V)$.  
\endproclaim

\proclaim{\smc Theorem 1.1} 
Let $p$ be a prime, let $V$ be a vector space of dimension $d$ over
$\fp$, and let $\K$ stand for the $n$th Morava K-theory.

\noindent
(a) If $p$ is odd, then $\K(BV)$ is not a graded permutation module
for $GL(V)$.

\noindent
(b) If $p=2$ and $n=1$, then for any $d$, 
$\K(BV)$ is a graded permutation module for $GL(V)$.

\noindent
(c) For $p=2$, $n>1$ and $d\geq 4$, $\K(BV)$ is 
not a graded permutation module for $GL(V)$ 
if $d$ is greater or equal to the smallest 
prime divisor of $n$.

\noindent
(d) For $p=2$ and $d=3$, $\K(BV)$ is not a graded permutation 
module if  $n$ is a multiple of three, or if $n$ is $2$, $4$, or $5$.

\noindent
(e) For $p=2$ and $d=2$, $\K(BV)$ is a graded permutation module 
for $GL(V)$ if and only if $n$ is odd.  

\endproclaim

\proclaim{\smc Theorem 1.2} The $\K[GL(V)]$-modules $\K(BV)$ and 
$\K[\Hom(V,(\fp)^n)]$ are (ungraded) isomorphic in the following
cases: 

\noindent 
(a) For $n=1$, for any $p$ and $d$.

\noindent 
(b) For $d=2$, $p=2$, and any $n$.  

\noindent
And are isomorphic as $\K[SL(V)]$-modules in the case:

\noindent 
(c) $d=2$, $p=3$, $n=1$, $2$ or $3$.  

\endproclaim 

\proclaim{\smc Theorem 1.3} $\K(BV)$ is not a permutation module for
$\K[U(V)]$ in the following cases:  

\noindent 
(a) $d=3$, $p=3$, $n=2$, 

\noindent 
(b) $d=3$, $p=5$, $n=2$.  

\noindent
In the following cases, as well as those implied by Theorems 1.1~and~1.2, 
$\K(BV)$ is a graded permutation module for $\K[U(V)]$: 

\noindent
(c) $d=3$, $p=2$, $n=2$, $3$ or $4$.  

\endproclaim 

Work of Kriz [\kri] shows that for any $V$, $\K(BV)$ is a (graded)
permutation module for any subgroup of $GL(V)$ of order $p$.  In the
cases covered by Theorem~1.3, the group $U(V)$ has order $p^3$, and for
$p>2$ it contains no element of order $p^2$.  The gap between
Theorem~1.3 and a special case of Kriz's result is filled by: 

\proclaim{\smc Theorem 1.4} Let $d=3$, let $p=3$ or $5$, and let $H$ be any
subgroup of $GL(V)$ of order $p^2$.  Then $K(2)^*(BV)$ is not a
permutation module for $H$.  
\endproclaim

Statements 1.1(b) and 1.2(a) are corollaries of Kuhn's description of the
mod-$p$ K-theory of $BG$ [\kui].  Our interest 
in these questions was aroused by [\bru], in which it is shown that 
in the case $V= (\zz/2)^2$, $\K(BV)$ is a graded permutation module
for $n=3$ but is not a graded permutation module for $n=2$, i.e., 
the cases $n=2$ and $n=3$ of 1.1(e).

The remaining sections of the paper are organised as follows.  
In Section~2 we describe $\K(BV)$, the action of $GL(V)$, 
and the process of reduction to questions concerning
finite-dimensional $\fp$-vector spaces.  This material is well-known
to many topologists, but we hope that its inclusion will make the 
rest of the paper accessible to the reader who knows nothing about 
Morava K-theory.  We also include some remarks concerning 
$\K[\Hom(V,(\fp)^n)]$.  In Section~3 we prove those of our results 
that require only Brauer character methods.  In Section~4 we deduce 
1.1(b)~and~1.2(a) from Kuhn's work on mod-$p$ K-theory and give a second 
proof of~1.1(b).  Section~5 studies $\ftwo[GL_2(\ftwo)]$-modules, 
and contains proofs of 1.1(e)~and~1.2(b).  In Section~6 
we describe how to decompose $\fthree[SL_2(\fthree)]$-modules, 
and outline the proof of~1.2(c).  In Section~7 we describe 
our algorithm for determining when a module for a $p$-group 
is a permutation module, and outline the proofs of the rest of the
results we have obtained using computer calculations.  The algorithm 
of Section~7 may be of independent interest.  Section~8 contains 
the tables of computer output relevant to Sections 6~and~7, together 
with some final remarks.

\noindent

\head 2. {\it Preliminaries} \endhead

Fix a prime $p$.  The $n$th Morava K-theory, $\K$ (which depends on
$p$ as well as on the positive integer $n$), is a generalized 
cohomology theory whose coefficient ring is the ring
$\fp[v_n,v_n^{-1}]$ of Laurent polynomials in $v_n$, which has degree 
\hbox{$-2(p^n-1)$}. All graded modules for this ring are free, which implies
that there is a good K\"unneth theorem for $\K$.  For any graded
$\K$-module $M$, let $\overline M$ be the quotient \hbox{$M/(1-v_n)M$}.  
Then $\overline M$ is an $\fp$-vector space, naturally graded by the
cyclic group \hbox{$\zz/2(p^n-1)$}.  If $M$ is a graded $\K[G]$-module for
some finite group $G$, then $\overline M$ is naturally a
$\zz/2(p^n-1)$-graded $\fp[G]$-module, and $M$ is determined up to
isomorphism by $\overline M$.  It is easy to see that $M$ is a
(graded) permutation module for $\K[G]$ if and only if $\overline M$
is a ($\zz/2(p^n-1)$-graded) permutation module for $\fp[G]$.  

For $C$ a cyclic group of order $p^m$, it may be shown [\rawi] 
that the Morava K-theory of $BC$ is a truncated polynomial ring on 
a generator of degree two: 
$$\K(BC)= \K[x]/(x^{p^{mn}}).$$
The generator $x$ is a Chern class in the sense that it is the image
of a certain element of $K(n)^2(BU(1))$ under the map
induced by an inclusion of $C$ in the unitary group $U(1)$.  From the 
K\"unneth theorem mentioned above it follows that if $V$ is an
elementary abelian $p$-group of rank $d$, then 
$$\K(BV)=\K[x_1,\ldots,x_d]/(x_1^{p^n},\ldots,x_d^{p^n}),$$
where $x_1,\ldots,x_d$ are Chern classes of $d$ 1-dimensional
representations of $V$ whose kernels intersect trivially.  The Chern
class of a representation is natural, and the $d$ representations
taken above must generate the representation ring of $V$.  Thus the
action of $GL(V)$ on $\K(BV)$ may be computed from its action on
$\Hom(V,U(1))$ together with an expression for the Chern class of a
tensor product $\rho\otimes\theta$ of two 1-dimensional
representations in terms of the Chern classes of $\rho$ and $\theta$.

For any generalized cohomology theory $h^*$ such that $h^*(BU(1))$ is  
a power series ring $h^*[[x]]$ (Morava K-theory has this property), 
Chern classes may be defined, and there is a power series $x\fs y
\in h^*[[x,y]]$ expressing the Chern class of a tensor product of line
bundles in terms of the two Chern classes.  This power series is
called the formal group law for $h^*$, because it satisfies the axioms
for a 1-dimensional commutative formal group law over the ring $h^*$.
Since each Chern class in $\K(BV)$ is nilpotent of class $p^n$, we
need only determine $x\fs y$ modulo $(x^{p^n}, y^{p^n})$.  This is the
content of the following proposition, which is well-known, but for
which we can find no reference.  

\proclaim{\smc Proposition 2.1}  
Modulo the ideal generated by $x^{p^n}$ and $y^{p^n}$, the formal sum 
$x\fs y$ for $\K$ is 
$$x\fs y = x + y - v_n(\sum_{i=1}^{p-1} 
{1\over p} {p\choose i} x^{ip^{n-1}} y^{(p-i)p^{n-1}}).$$
\endproclaim 

\demo{Sketch proof}
First we recall the formal sum for $BP^*$, Brown-Peterson cohomology
[\wil].  Let $l$ be the power series 
$$l(x) = \sum_{i\geq 0} m_ix^{p^i},$$
where $m_0=1$, but the remaining $m_i$'s are viewed as indeterminates,
and let $e(x)$ be the compositional inverse to $l$, i.e., a power
series such that $e(l(x))=l(e(x))=x$.  The $BP^*$ formal sum is the
power series $e(l(x)+l(y))$.  The $\K$ formal sum may be obtained as
follows:  Take the $BP^*$ formal sum, replace the indeterminates $m_i$
by indeterminates $v_i$ using the relation 
$$v_j = pm_j - \sum_{i=1}^{j-1}m_iv_{j-i}^{p^i},$$
set $v_i=0$ for $i\neq n$, by which point all the coefficients lie in
$\zz_{(p)}$, and take the reduction modulo $p$.  To calculate the $\K$
formal sum, it is helpful to set $v_i=0$ for $i\neq n$ as early as
possible, and one may as well set $v_n=1$, since every term in $x\fs
y$ has degree 2.  Solving for the $m_i$'s in terms of the $v_i$'s
gives 
$$\eqalign{m_i &= 0\quad\hbox{if $n$ does not divide $i$,}\cr
m_{ni} &= 1/p^i.}$$
Thus to compute $x\fs y$, let $e'(x)$ be the compositional inverse to 
$$l'(x)=\sum_{i\geq 0} x^{p^{ni}}/p^i,$$
and then $x\fs y$ is the mod-$p$ reduction of $e'(l'(x)+l'(y))$.  It
is easy to see that 
$$e'(x)\equiv x-x^{p^n}/p \quad\hbox{modulo $x^{2p^n}$,}$$
and so 
$$x\fs y \equiv x + y - (x + y)^{p^n}/p \quad\hbox{modulo
$(x^{p^n},y^{p^n}, p)$}.$$
The claimed result follows.  %\qed
\enddemo

Using the reduction $M \mapsto \overline M$ as at the start of this
section and Proposition~2.1, the study of the graded $\K[GL(V)]$-module
structure of $\K(BV)$ reduces to the study of the $\zz/(p^n-1)$-graded
$\fp[GL_d(\fp)]$-module $K_{n,d}^*$ defined below.  

As an $\fp$-algebra, 
$$K_{n,d}^*\cong \fp[x_1,\ldots,x_d]/(x_1^{p^n},\ldots,x_d^{p^n}).$$
Each $x_i$ has degree 1, and the $GL_d(\fp)$-action is compatible
with the product.  
The action of the matrix $(a_{ij})\in GL_d(\fp)$ is given by 
$$x_j\mapsto e'\bigl(\sum_i l'(x_i)a_{ij}\bigr),$$
where $e'$ and $l'$ are as in the proof of Proposition~2.1.  (Recall that
for any $V$ we take $GL(V)$ to act on the right of $V$, and hence
obtain a left $\K[GL(V)]$-module structure on $\K(BV)$.)  

Note that we have halved the original degrees because $\K(BV)$ is
concentrated in even degrees.  Until recently it was an open problem
whether a similar statement holds for arbitrary finite groups,
although some cases had been verified [\rawi,\jrh,\teyai,\teyaii,\bjoe]. 
Kriz has recently announced that this is not the case [\kri].  

If we are only interested in Brauer characters, or equivalently
composition factors, then a further simplification may be made, see
[\kuh].   Let $L_{n,d}^*$ denote the algebra of polynomial functions 
on $(\fp)^d$, modulo the ideal of $p^n$th powers of elements of
positive degree.  Grade $L_{n,d}^*$ by $\zz/(p^n-1)$, and let
$GL_d(\fp)$ act on $L_{n,d}^*$ by its natural action on the
polynomial functions.  Thus $L_{n,d}^*$ is a truncated polynomial
algebra $\fp[x_1,\ldots,x_d]/(x_1^{p^n},\ldots,x_d^{p^n})$, cyclically
graded, and having the standard action of $GL_d(\fp)$.  

\proclaim{\smc Lemma 2.2} $K^*_{n,d}$ has a series of (graded) submodules such
that the direct sum of the corresponding quotients is isomorphic to
$L^*_{n,d}$.  In particular, $K^*_{n,d}$ and $L^*_{n,d}$ have the same 
composition factors (as graded modules).
\endproclaim

\demo{Proof} For each degree $k$, take the basis consisting of monomials
of length congruent to $k$ modulo $p^n-1$, and arrange them in blocks
with respect to length. For any $g$, the matrix of its action on $L^k_{n,d}$
with respect to this basis consists of square blocks along the diagonal,
whereas the corresponding matrix for the action on $K^k_{n,d}$ has some extra 
entries below the blocks.  %\qed
\enddemo

The permutation module $\K[\Hom(V,(\fp)^n)]$ occurs in the statement
of Theorem~1.2, so we complete this preliminary section with some 
remarks concerning this module.  
If $\phi$ is a homomorphism from $V$ to 
$\fp^n$, then $g\in GL(V)$ acts by composition, i.e., 
%$$\phi^g(v)= \phi(gv).$$
$$g\phi(v) = \phi(vg).$$
Since we view $GL(V)$ as acting on the right of $V$, this makes
$\Hom(V,(\fp)^n)$ into a left $GL(V)$-set.  The 
$GL(V)$-orbits in $\Hom(V,(\fp)^n)$ may be described as follows.  
For $W$ a subspace of $V$, let $H(W)\leq GL(V)$ be 
$$H(W) = \{ g\in GL(V) :  vg - v \in W \,\, \forall v\in V\}.$$
For example, $H(\{0\}) = \{1\}$, and $H(V)=GL(V)$.  
For $0\leq i \leq \dim(V)$, let $H_i$ be $H(W_i)$ for some $W_i$ of 
dimension $i$.  Thus $H_i$ is defined only up to conjugacy, but this 
suffices to determine the isomorphism type of the $GL(V)$-set 
$GL(V)/H_i$.  Now let $\phi$ be an element of $\Hom(V,(\fp)^n)$.  
The stabilizer of $\phi$ in $GL(V)$ is the subgroup $H(\ker(\phi))$, 
and the orbit of $\phi$ consists of all $\phi'$ such that 
$\Im(\phi')=\Im(\phi)$.  
It follows that as $GL(V)$-sets,  
$$\Hom(V,(\fp)^n)\cong \coprod_{0\leq i\leq\dim(V)} m(n,i)\cdot G/H_i,$$
where $m(n,i)$ is the number of subspaces of $(\fp)^n$ of dimension $i$.
Thus to decompose the module $\fp[\Hom(V,(\fp)^n)]$, it suffices to 
decompose each $\fp[GL(V)/H_i]$.

\head 3. {\it Brauer characters} \endhead

In this section we shall prove most of the negative results of
Theorem~1.1.  Firstly, we describe how to compute the values of the modular 
characters afforded by the modules $L^k_{n,d}$.  
As a general reference, see [\cur], in particular \S 17.
Fix an embedding of the multiplicative group of the algebraic closure
of $\fp$ in the group of roots of 1 in $\cc$.  
Let $g$ be a $p$-regular element of $GL(V)$, 
i.e., an element whose order is coprime to $p$, 
and let $\lambda_1, \lambda_2, 
\ldots ,\lambda_d$ denote the images in $\cc$ of 
the eigenvalues of its action on $V^*$. Then the
Brauer character of $g$ is
$$
\chi_{V^*}(g)= \lambda_1+\lambda_2+\cdots+\lambda_d\,.
$$
Two $\fp [G]$-modules have the same Brauer character if and only if they
have the same composition factors.  
To compute the character of $L_{n,d}^k$ we proceed as follows: an argument 
similar to the one used to prove Molien's theorem (see e.g. [\cur], p.~329) 
shows that the character of a truncated polynomial algebra has a generating 
function
$$
  f_{g}(t)=\prod_{i=1}^{d}
   \Bigl(\frac{1-(\lambda_it)^{p^n}}{1-\lambda_it}\Bigr)\,.
$$
Then the character of $L_{n,d}^*$ evaluated at $g$ is simply $f_{g}(1)$, 
whereas for each degree $k$ (recall that we are grading cyclically) one has
$$
\chi_{L_{n,d}^k}(g)=\frac{1}{p^n-1}\sum_\tau \tau^{-k}f_{g}(\tau)\,,\tag3.1
$$
where the sum ranges over all $(p^n-1)$-st roots of unity---to see
this, recall that the sum, over all $m$th roots of unity $\lambda$, of 
$\lambda^k$ is equal to zero if $m$ does not divide $k$, and equal to 
$m$ if $m$ does divide $k$.  

\demo{Proof of 1.1(a)} 
Let $D$ be the subgroup of diagonal matrices in $GL_d(\fp)$, so that 
$D$ is isomorphic to a direct product of $d$ cyclic groups of order
$p-1$.  In $K^*_{n,d}$, each monomial in $x_1,\ldots,x_d$ is an
eigenvector for $D$, and the monomials fixed by $D$ are those in which
the exponent of each $x_i$ is divisible by $p-1$.  Hence if $p-1$ does
not divide $k$, then $K^k_{n,d}$ cannot be a permutation module for
$D$ because it contains no $D$-fixed point.  %\qed
\enddemo

\demo{Proof of 1.1(c)}
As already said above, this is done by computing the character values on
certain 2-regular elements of $GL(V)$. 
We shall first look at the case where $d$ equals a prime divisor $q$ of the
fixed number 
$n$. Consider an element, $g_q$ say, of $GL_q(\ftwo)$ which permutes the
$2^q-1$ nontrivial elements of $(\ftwo )^q$ cyclically. (To see that there is
always such an element consider the action of the multiplicative group of
${\bold F}_{2^q}$ on the additive group of ${\bold F}_{2^q}$.)  
The set of eigenvalues of $g_q$ contains a primitive $(2^q-1)$st root
of unity, and is closed under the action of the Galois group
$\Gal({\bold F}_{2^q}/\ftwo )$.  Hence the Brauer lifts of the
eigenvalues of $g_q$ are $\lambda$,
$\lambda^2,\ldots,\lambda^{2^{q-1}}$
for some primitive $(2^q-1)$-st root of unity $\lambda\in \cc$. 
Consequently, the generating function for the character afforded 
by $L_{n,q}^*$ is given by
$$f_{g_q}(t)=\prod_{i=0}^{q-1}
  \Bigl(\frac{1-(\lambda^{2^i}t)^{2^n}}{1-\lambda^{2^i}t}\Bigr)\,.$$
If $\tau$ is a $(2^n-1)$-st root of unity, one gets 
$$
f_{g_q}(\tau)=
\left\{
 \aligned
   2^n\phantom{1}\quad &\quad
   \text{ if $\tau\in\{\lambda^{-2^i},\quad i=0,1,\ldots q-1\}$}\\
   1\phantom{2^n}\quad&\quad\text{ otherwise.}
  \endaligned
\right.$$
Thus evaluating the formula (3.1) for the character 
afforded by $L_{n,q}^k$ yields
$$\chi_{L_{n,q}^k}(g_q)=\frac{1}{2^n-1}\sum_{\tau\ne\lambda^{-2^i}}\tau^k
   +\frac{2^n}{2^n-1}\sum_{i=0}^{q-1}\lambda^{2^ik}$$   
which is equal to 
$$\left\{
  \aligned
    q+1\quad&\quad\text{ for $k=0$}\\
    \sum_{i=0}^{q-1}\lambda^{2^ik} &\quad\text{ for $k\ne 0$.}
   \endaligned
\right.$$

\noindent
Specializing to the case $k=1$, this sum is never zero, since the 
powers $\lambda^i$ for $i$ coprime to $2^q-1$ form a $\qq$-basis for 
$\qq[\lambda]$.  For $q>2$ the sum is not a rational, because it is
not fixed by the whole Galois group $\Gal(\qq[\lambda]/\qq)$.  In the
case $q=2$, one obtains $-1$ (the sum of the two primitive third roots
of unity).  
Since permutation modules have positive integer character values, this shows
that $K(n)^*(BV_q)$ is not a graded $GL_q(\ftwo )$-permutation module. 
To proceed with vector spaces of dimension bigger than $q$ we consider the
cases $q>2$ and $q=2$ separately. In the first case
we use the following lemma to conclude that the
character still takes non-integer values on certain elements of $GL(V)$.
Let $g$ be an (arbitrary) 2-regular element of $GL_d(\ftwo )$ and $I_r$ the 
$r\times r$ identity matrix.
If we denote by $g\times I_r$ the element of $GL_{d+r}(\ftwo )$ 
which acts like $g$ on the first $d$ generators of $L_{n,d+r}^*$ 
and trivially on the last $r$, one has
\proclaim{\smc Lemma 3.2} 
$\chi_{L_{n,d+r}^k}(g\times I_r)=\chi_{L_{n,d}^k}(g)+
   %\bigl(\frac{2^{nr}-1}{2^n-1}\bigr)\chi_L(g)$ .
    \bigl((2^{nr}-1)/(2^n-1)\bigr)\chi_L(g)$ .
\endproclaim
\demo{Proof}
The generating function for $g\times I_r$ is obtained from the one for $g$
as the product with $r$ factors $(1+t+t^2+\ldots +t^{2^n-1})$, thus
$$
\aligned
   \chi_{L_{n,d+r}^k}(g\times I_r)&=
            \frac{1}{2^n-1}\sum_\tau \tau^{-k} f_{g\times I_r} (\tau)\\
   &=\frac{1}{2^n-1}\sum_{\tau\ne 1}\tau^{-k}f_{g}(\tau)
      \Bigl(\frac{1-\tau^{2^n}}{1-\tau}\Bigr)^r
          +\frac{2^{nr}}{2^n-1}f_{g}(1)\\
   &= \chi_{L_{n,d}^k}(g)+
   \bigl(\frac{2^{nr}-1}{2^n-1}\bigr)f_{g}(1)\, .%\qed
\endaligned
$$
\enddemo
Thus $g_q\times I_r$ will do the trick when $V$ has rank $q+r$. This
fails for $q=2$, whence we choose the element $g^\prime$ which consists of 
$d/2$ copies of $g_2$ = $0\,1\choose1\,1$ arranged along the diagonal if $d$ 
is even, and add
an extra diagonal entry 1 if $d$ is odd. Then a computation similar to the 
one carried out in the previous lemma shows that for $k\not\equiv 0$ mod 3,
$$\chi_{L_{n,d}^k}(g^\prime)=
   \left\{
     \aligned
	-\frac{2^{nd/2}-1}{2^n-1}&\quad\text{ if $d$ is even}\\
	-\frac{2^{n(d-1)/2}-1}{2^n-1}+1&\quad\text{ if $d$ is odd}
    \endaligned
   \right.
$$
For $d>3$ these numbers are negative.  
%\qed
\enddemo

The other parts of Theorem~1.1 that may be proved using Brauer
characters are some cases of 1.1(d) and both implications of 1.1(e).  
The details are similar to the above proof so we shall not give them.  
In the case when $V$ has rank 3, evaluation of the Brauer character of an
element of $GL(V)$ of order 7 shows that $L^1_{n,3}$ is not a
$GL(V)$-permutation module if 3 divides $n$. 
Similarly, when $V$ has rank 2, the Brauer character of an element of
$GL(V)$ of order 3 on $L^1_{n,2}$ is negative if $n$ is even.  
When $V$ has rank~2 and $n$ is odd, it may be shown that for each $k$,
any $GL(V)$-module having the same Brauer character as $L^k_{n,2}$ is
a permutation module.  This shows that for $n$ odd, $K^k_{n,2}$ is a
permutation module, but does not specify which one. In Section~5 we
shall describe the isomorphism type of $K^k_{n,2}$ and $L^k_{n,2}$ for
all $n$~and~$k$, giving an alternative proof of 1.1(e).

\head 4. {\it On $K(1)$ }\endhead
Here we describe how those parts of Theorems 1.1~and~1.2 that concern
$K(1)^*$ (i.e., 1.1(b) and 1.2(a)) follow from Kuhn's description of the
mod-$p$ K-theory of finite groups [\kui].  An \lq elementary' proof, 
working directly with the description of $K^*_{1,d}$ in the previous
section, would be more in keeping with the rest of the paper.  We give
such a proof in the case $p=2$.  

\demo{Proof of 1.1(b) and 1.2(a)} First, note that for $p=2$, $v_1$ has
degree $-2$, so that the \lq cyclically graded' modules $K^*_{1,d}$
are in fact concentrated in a single degree.  Hence 1.2(a) implies 1.1(b).
To prove 1.2(a) recall [\wil] that the spectrum representing mod-$p$
K-theory splits as a wedge of one copy of each of the 0th, 
$2\hbox{nd},\ldots,(2p-4)\hbox{th}$ suspensions of the spectrum
representing $K(1)^*$.  Since $K(1)^*(BV)$ is concentrated in  
even degrees it follows that
$K^*_{1,d}$ is naturally isomorphic to $K^0(BV;\fp)$.  In [\kui] it is
shown that for any $p$-group $G$, $K^0(BG;\fp)$ is naturally
isomorphic to $\fp\otimes R(G)$, where $R(G)$ is the (complex)
representation ring of $G$.  The case $G=V$ gives 1.2(a), because as a 
$GL(V)$-module, $\fp\otimes R(V)$ is isomorphic to $\fp[\Hom(V,\fp)]$.
%\qed
\enddemo

\demo{Alternative proof, $p=2$} In this case, $K^*_{1,d}$ is
isomorphic to an exterior algebra $\Lambda[x_1,\ldots,x_d]= 
\ftwo [x_1,\ldots,x_d]/(x_i^2)$.  The monomial 1 generates a trivial
$GL(V)$-summand.  Let $H$ be the subgroup of $GL(V)$ fixing $x_1$.
Then $H$ is the subgroup of $GL(V)$ stabilizing some hyperplane $W$  
and inducing the identity map on the quotient $V/W$.  
There is a $GL(V)$-set isomorphism 
$$\Hom(V,\ftwo )\cong GL(V)/GL(V) \amalg GL(V)/H,$$
so it will suffice to show that the submodule $M$
generated by $x_1$ contains each monomial in $\Lambda[x_1,\ldots,x_d]$
of strictly positive length.  The permutation matrices permute the
monomials of any given length transitively.  Assume that $M$ contains
all the monomials of length $i$ (this holds for $i=1$), and let $g\in
GL(V)$ be such that 
$$gx_1= x_1,\ldots,gx_{i-1}=x_{i-1},\quad gx_i = x_i\fs x_{i+1}.$$
Then 
$$g(x_1\ldots x_i) + x_1\ldots x_i + x_1\ldots x_{i-1}x_{i+1} 
= x_1\ldots x_ix_{i+1} \in M,$$
so $M$ contains all monomials of length $i+1$.  
%\qed
\enddemo

It should be possible to give an \lq elementary' proof of 1.2(a) for
$p>2$ by considering the element $x_1+ x_1^2+\cdots + x_1^{p-1}$, but
we have not done so.

\head 5. {\it  When $V$ has order four }\endhead

%In this section
Here we shall prove Theorems 1.1(e) and 1.2(b), which concern 
$\K(BV)$ for $V$ of dimension two over $\ftwo $.  We determine the
structure of $L^k_{n,2}$ as a $GL_2(\ftwo )$-module, and deduce that 
$K^k_{n,2}$ and $L^k_{n,2}$ are isomorphic.  Note that it is also possible
to prove 1.1(e) using the methods of Section~3 without determining the 
isomorphism type of $K^k_{n,2}$.  Throughout this section, let
$V=(\ftwo)^2$.  

There are three isomorphism types of indecomposable $\ftwo 
[GL(V)]$-modules:  the 1-dimensional trivial module $T$; the natural
module $V$ which is both simple and projective (and is the Steinberg
module for $GL(V)$); and a module $N$ expressible as a non-split
extension of $T$ by $T$, which is the projective cover of $T$.  
Each module is self-dual.  There are four conjugacy classes of
subgroups of $GL(V)$.  The transitive permutation modules are the
following four modules:  
$$T,\qquad N,\qquad T\oplus V, \qquad N\oplus 2V.$$
Let $S^*[V^*]$ stand for the algebra of polynomial functions on $V$ as
a graded $GL(V)$-module.  

\proclaim{\smc Proposition 5.1} The generating functions $P_T$, $P_N$ and
$P_V$ for the number of each indecomposable $GL(V)$-summand of
$S^*[V^*]$ are the following power series:  
$$P_T(t)={1\over 1-t^2},\qquad 
P_N(t) = {t^3 \over (1-t^2)(1-t^3)},\qquad
P_V(t) = {t\over (1-t)(1-t^3)}.$$
\endproclaim

\demo{Proof} 
Recall that the ring of invariants $S^*[V^*]^{GL(V)}$ is a free
polynomial ring on two generators of degrees two and three (see
[\wilk]).  The Poincar\'e series for $S^*[V^*]$, and the ring of
invariants, together with the generating function for the 
Brauer character of an element of $GL(V)$ of order three give the
three equations below, whose solution is as claimed.  
$$\eqalign{
P_T + 2P_N + 2P_V &= {1\over (1-t)^2}\cr
P_T + P_N &= {1\over (1-t^2)(1-t^3)}\cr
P_T + 2P_N - P_V &= {1-t \over 1-t^3}
\hbox to 0pt{\hbox to 4cm{\hfill\phantom{a}%\qed
}}\cr}$$
\enddemo

\proclaim{\smc Proposition 5.2} Let $k$ be an element of $\zz/(2^n-1)$.  
The direct sum decomposition for the module $L^k_{n,2}$ is:  
%$$
%\aligned
%2T\oplus {2^n-2\over 6}N \oplus {2^n+1\over 3}V 
%&\quad\hbox{for $n$ odd, $k=0$,}\\
%T\oplus {2^n-2\over 6}N \oplus {2^n+1\over 3}V 
%&\quad\hbox{for $n$ odd, $k\neq 0$,}\\
%2T\oplus {2^n+2\over 6}N \oplus {2^n-1\over 3}V 
%&\quad\hbox{for $n$ even, $k=0$,}\\
%T\oplus {2^n+2\over 6}N \oplus {2^n-1\over 3}V 
%&\quad\hbox{for $n$ even, $k\neq0$, $k\equiv 0$ mod 3,}\\
%T\oplus {2^n-4\over 6}N \oplus {2^n+2\over 3}V  
%&\quad\hbox{for $n$ even, $k\not\equiv0$ mod 3.}\\
%\endaligned
%$$
$$
\aligned
2T\oplus (2^n-2)/ 6\, N \oplus (2^n+1)/ 3\, V 
&\quad\hbox{for $n$ odd, $k=0$,}\\
T\oplus (2^n-2)/ 6\,N \oplus (2^n+1)/ 3\,V 
&\quad\hbox{for $n$ odd, $k\neq 0$,}\\
2T\oplus (2^n+2)/ 6\,N \oplus (2^n-1)/ 3\,V 
&\quad\hbox{for $n$ even, $k=0$,}\\
T\oplus (2^n+2)/ 6\,N \oplus (2^n-1)/ 3\,V 
&\quad\hbox{for $n$ even, $k\neq0$, $k\equiv 0$ mod 3,}\\
T\oplus (2^n-4)/ 6\,N \oplus (2^n+2)/ 3\,V  
&\quad\hbox{for $n$ even, $k\not\equiv0$ mod 3.}\\
\endaligned
$$
\endproclaim

\demo{Proof} Let $\tilde L^*$ be the truncated symmetric algebra
$L^*_{n,2}$, but graded over the integers rather than over the
integers modulo $2^n-1$. Then $\tilde L^k=\{0\}$ for $k> 2(2^n-1)$,
and for $k= 2(2^n-1)$, $\tilde L^k$ is isomorphic to $T$, generated by
$x_1^{2^n-1}x_2^{2^n-1}$.  For $0<k<2^n-1$, viewing $k$ as either an
integer or an integer modulo $2^n-1$ as appropriate, $L^k_{n,2}$ is
isomorphic to $\tilde L^k \oplus \tilde L^{2^n-1 + k}$, while 
$L^0_{n,2}$ is isomorphic to $\tilde L^0\oplus \tilde L^{2^n-1}
\oplus \tilde L^{2(2^n-1)}\cong 2T\oplus \tilde L^{2^n-1}$.  For 
$0\leq k\leq 2^n-1$, $\tilde L^k$ is isomorphic to $S^k[V^*]$.  The
product structure on $\tilde L^*$ gives a duality pairing
$$\tilde L^k \times \tilde L^{2(2^n-1)-k}\rightarrow 
\tilde L^{2(2^n-1)}\cong T,$$
and since all $GL(V)$-modules are self-dual it follows that for
\hbox{$2^n\leq k \leq 2(2^n-1)$}, $\tilde L^k\kern 1pt\cong\kern 1pt 
S^{2(2^n-1)-k}[V^*]$.  
The claimed description of $L^k_{n,2}$ follows from 
Proposition~5.1.\kern 5pt%\qed
\enddemo

\proclaim{\smc Corollary 5.3} For each $n$ and each $k\in \zz/(2^n-1)$, 
$K^k_{n,2}$ and $L^k_{n,2}$ are isomorphic.  
\endproclaim 

\demo{Proof} For $k\neq 0$ $K^k_{n,2}$ has odd dimension, so must
contain at least one direct summand isomorphic to $T$.  It is easy to
see that $1$ generates a summand of $K^0_{n,2}$, and the same
dimension argument applied to a complement of this summand shows that
$K^0_{n,2}$ contains at least two summands isomorphic to $T$.  On the
other hand, $V$ and $N$ are projective, and %(Proposition~2.1)
(Lemma~2.2)
$K^k_{n,2}$ has a filtration such that the sum of the factors is
isomorphic to $L^k_{n,2}$.  Hence $K^k_{n,2}$ has at least as many 
$N$ summands and $V$ summands as $L^k_{n,2}$.  This accounts for all
the summands of $K^k_{n,2}$.  
%\qed
\enddemo

The proof of 1.1(e) follows easily from the given description
of $K^*_{n,2}$.  For 1.2(b), recall from the end of Section~2 that 
%$$\eqalign{
%\Hom(V,(\ftwo)^n)= GL(V)/GL(V) &\amalg 
%(2^n-1)GL(V)/H_1 \amalg \cr
%&(2^n-1)(2^n-2)/6 GL(V)/\{1\},\cr}$$
$$\eqalign{
  \Hom(V,(\ftwo)^n)= &\,GL(V)/GL(V) \cr
  &\quad\amalg\, (2^n-1)\cdot GL(V)/H_1\cr
  &\quad\amalg\, (2^n-1)(2^n-2)/6\cdot GL(V)/\{1\},\cr}$$
where $H_1$ is a subgroup of $GL(V)$ of order two, 
and that 
$$\ftwo[GL(V)/H_1]\cong T\oplus V,\qquad 
\ftwo[GL(V)/\{1\}]\cong N\oplus 2V.$$

The argument used in the proof of Corollary~5.3 shows that for any $p$,
$n$, and $k$ such that $L^k_{n,2}$ contains at most one non-projective
summand, $K^k_{n,2} \cong L^k_{n,2}$.  For odd primes this does not
always occur however.  If $0<k<p^n-1$ then $L^k_{n,2}$ splits as a
direct sum of submodules of dimensions $k+1$ and $p^n-k-2$ coming from
the standard $\zz$-grading on the truncated polynomial algebra.  If $k$
is not congruent to either $-1$ or $-2$ modulo $p$, the dimensions of
these summands are not divisible by $p$, and hence $L^k_{n,2}$ contains
at least two non-projective indecomposable summands.  The calculations
described in the next section show that for $p=3$, $n=2,3$ and
$0<k<3^n-1$, the module $K^k_{n,2}$ has exactly one non-projective
indecomposable summand.  It follows that for $p=3$, $K^k_{n,2}$ and
$L^k_{n,2}$ are not necessarily isomorphic.  

\head 6. {\it When $V$ has order nine }\endhead

Our results concerning the $SL_2(\fthree )$-module structure of $K^*_{n,2}$
in the case when $p=3$ were obtained by computer.  We wrote a Maple program
to generate matrices representing the action of a pair of generators
for $SL_2(\fthree )$ on $K^k_{n,2}$.  These matrices were fed to a GAP
[\gap] program which, given a matrix representation of $SL_2(\fthree )$,
outputs a list of its indecomposable summands.  In fact the output
from the Maple program needed a little editing before being read into 
the GAP program.  This was done by a third program, although it could
equally have been done by hand.  

Using standard techniques of representation theory [\alp,\cur], the
following facts may be verified.  For $V= (\fthree)^2$, there are three
simple $\fthree[SL(V)]$-modules: the trivial module $T$, the natural
module $V$, and a simple projective module $P=S^2(V)$ of dimension
three.  There are three blocks.  The blocks containing $T$ and $V$
each contain three indecomposable modules, each of which is uniserial.
This data may be summarised as follows: 
$$\eqalign{
\hbox{block of } T=I_1:&\qquad T\rightarrowtail I_2\twoheadrightarrow
T,\qquad T\rightarrowtail I_3 \twoheadrightarrow I_2,\cr
\hbox{block of } V=I_4:&\qquad V\rightarrowtail I_5\twoheadrightarrow
V,\qquad V\rightarrowtail I_6 \twoheadrightarrow I_5,\cr
\hbox{block of $P= I_7$}:&\qquad 
\hbox{contains no other indecomposables.}\cr}$$

Letting $\tau$ stand for the element of order two in $SL(V)$ and 
$\sigma$ for the sum of the six elements of $SL(V)$ of order four, the
block idempotents are 
$$b_T = 2 + 2\tau + 2\sigma,\qquad 
b_V = 2 + \tau,\qquad 
b_P = \sigma.$$
The modules in any single block are distinguishable by their
restrictions to a cyclic subgroup of $SL(V)$ of order three.  Thus if
$\alpha$ is an element of $SL(V)$ of order three, and $M$ is an
$SL(V)$-module, the direct summands of $M$ are determined by the ranks
of the elements of $\End(M)$ representing the actions of the following
seven elements of $\fthree[SL(V)]$:  
$$b_T,\quad (1-\alpha)b_T, \quad (1-\alpha)^2b_T, 
\quad b_V, \quad (1-\alpha)b_V, \quad (1-\alpha)^2b_V, \quad b_P.$$
More precisely, if the seven ranks are $r_1,\ldots,r_7$, and $n_i$
stands for the number of factors of $M$ isomorphic to $I_i$, then 
$$n_1 = r_1 - 2r_2 + r_3,\quad 
n_2 = r_2 - 2r_3, \quad
n_3 = r_3, \quad
n_4 = r_5 - 2r_6,$$
$$
n_5 = r_4 - 2r_5 + r_6,\quad
n_6 = (2r_5 - r_4)/2,\quad
n_7 = r_7/3.$$

Our GAP program reads in matrices representing the action on $M$ of
a certain pair of generators for $SL(V)$, 
and calculates $n_1,\ldots,n_7$ by first finding
$r_1,\ldots,r_7$ as above.

Recall from the end of Section~2 that there is an 
isomorphism of (right) $GL(V)$-sets:  
$$\eqalign{
\Hom(V,(\fthree)^n)=&GL(V)/GL(V) \cr
&\quad\amalg 
(3^n-1)/2\cdot GL(V)/H_1\cr 
&\quad\amalg (3^n-1)(3^n-3)/48\cdot GL(V)/\{1\},\cr}$$ 
where $H_1$ is the subgroup stabilizing a line $L$ in $V$ and acting
trivially on $V/L$.  As $SL(V)$-modules, it may be checked that 
$$\eqalign{
\fthree[GL(V)/GL(V)] &\cong I_1,\cr
\fthree[GL(V)/L(V)] &\cong I_1\oplus I_5\oplus I_7,\cr
\fthree[GL(V)/\{1\}] &\cong 2I_3\oplus 4I_6 \oplus 6I_7.
\cr}$$ From this information together with the results given in Table~8.1
it is easy to check the claim of Theorem~1.2(c).

There are fourteen indecomposable $GL(V)$-modules in four blocks, two
of which contain a single simple projective module.  The six
indecomposables in the block containing $V$ are comparatively hard to
distinguish, which is the reason why we considered only $SL(V)$.

\head 7. {\it Permutation modules for $p$-groups }\endhead
In this section we shall describe the computer programs used in 
the proofs of Theorem~1.3, Theorem~1.4, 
and the cases $n = 2$, 4, and 5 of 1.1(d).  
In Sections 5~and~6 our programs made use of the fact
that there were only finitely many indecomposable modules.  If $G$ is
a group whose Sylow $p$-subgroup is not cyclic, then $\fp[G]$ has
infinitely many indecomposable modules, so the same sort of methods
cannot work.  Here we shall describe an algorithm which may be used to
determine, for any $p$-group $G$, whether an $\fp[G]$-module is a
permutation module, and if so to decompose it.  (For a precise
statement, see Proposition~7.1 below.)  As before, we use a Maple
program to generate matrices representing the action of a
Sylow $p$-subgroup of $GL_d(\fp)$ on $K^k_{n,d}$ for various $p$, $k$,
$n$, and $d$, and we use a GAP program working with our algorithm to
decompose these modules.  

Our algorithm relies on the following fact [\cur]:  For $G$ a
$p$-group, any transitive permutation module for $\fp[G]$ has a
unique minimal submodule, which is the trivial module generated by the
sum of the elements of a permutation basis.  This implies that any
transitive permutation module is indecomposable.  Note that the
Krull-Schmidt theorem and the indecomposability of transitive
permutation modules together imply that if a graded $\fp[G]$-module
is a permutation module, then it is also a graded permutation module.

\proclaim{\smc Proposition 7.1} Let $G_1,\ldots,G_n$ be subgroups of a
$p$-group $G$, where the order of $G_{i+1}$ is at least the order of
$G_i$, and let $M$ be a (finitely generated) $\fp[G]$-module.  Let
$m_1,\ldots,m_n$ be the integers whose calculation is described below.
Then $M$ contains a submodule $M'$, where 
$$M'\cong m_1\fp[G/G_1]\oplus\cdots\oplus m_n\fp[G/G_n],$$ 
and $M'$ has maximal dimension among all submodules of $M$ isomorphic
to a direct sum of copies of the $\fp[G/G_i]$.  

To compute $m_i$, proceed as follows.   Let $M_0$ be the zero
submodule of $M$.  If $M_{i-1}$ has been defined, let 
$$M_i = M_{i-1} + 
\Im\Bigl(\bigl(\sum_{g\in G/G_i}g\bigr) : M^{G_i} \rightarrow M\Bigr),$$ 
where the sum ranges over a transversal to $G_i$ in $G$, $M^{G_i}$
denotes the $G_i$-fixed points of $M$, and the sum is an element of
$\fp[G]$ viewed as an element of $\End(M)$.  Now define 
$$m_i= \dim M_i - \dim M_{i-1}.$$

Without loss of generality, it may be assumed that no two of
$G_1,\ldots,G_n$ are conjugate.  The dimension of $M'$ is equal to the
sum $\sum_i m_i|G:G_i|$.  If $G_1,\ldots,G_n$ contains a
representative of each conjugacy class of subgroups of $G$, then $\dim
M'= \dim M$ if and only if $M$ is a permutation module.  
\endproclaim

\demo{Proof} First, recall that the socle, $\Soc(N)$, of a module $N$ is
the smallest submodule of $N$ containing every minimal submodule.  The
following statement is easy to prove, and will be useful below.  If
$L$ is a submodule of $M$, and $f:N\rightarrow M$ is a module
homomorphism, then $f$ is injective if and only if its restriction to
$\Soc(N)$ is injective.  If $f$ is injective, then $\Soc(f(N))=
f(\Soc(N))$, and the sum $L+f(N)$ in $M$ is direct if and only if 
the sum $\Soc(L)+f(\Soc(N))$ is direct.  

Module homomorphisms from $\fp[G/G_i]$ to $M$ are naturally bijective
with elements of $M^{G_i}$, where the element $x$ corresponds to the
homomorphism $\theta_x$ sending $1\cdot G_i$ to $x$.  The socle of
$\fp[G/G_i]$ is a trivial submodule generated by $\sum_{g\in G/G_i}
g\cdot G_i$, so its image under $\theta_x$ is generated by $\sum_{g\in
G/G_i} g\cdot x$.  It follows that any submodule of $M$ isomorphic to a
direct sum of copies of the  modules $\fp[G/G_1],\ldots,\fp[G/G_i]$ 
has socle contained in $M_i$, and in particular consists of at most 
$\dim M_i$ summands.  This shows that any submodule of $M$ isomorphic 
to a direct sum of $\fp[G/G_i]$'s has dimension less than or equal to 
$\sum_im_i|G:G_i|$, but it remains to exhibit a submodule $M'$ having
this dimension.  

Define $M'_0 $ to be the zero submodule of $M$, and assume that for
some $j$ with $1\leq j\leq n$ we have constructed a submodule
$M'_{j-1}$ of $M$ with 
$$M'_{j-1} \cong m_1\fp[G/G_1] \oplus \cdots \oplus
m_{j-1}\fp[G/G_{j-1}].$$
Let $x_1,\ldots,x_{m_j}\in M^{G_j}$ be such that the images
$\sum_{g\in G/G_j}g\cdot x_i$ form a basis for a complement to $M_{j-1}$ in
$M_j$.  Taking $L = M_{j-1}$, $N = m_j\fp[G/G_j]$, and $f:N\rightarrow
M$ the map sending the elements $(0,\ldots,1\cdot G_j,\ldots,0)$ to the
$x_i$'s, the statements in the first paragraph of the proof show that
$f$ is injective, and that $M'_j$ defined as the submodule of $M$
spanned by $M_{j-1}$ and the $x_i$'s is isomorphic to $M'_{j-1}\oplus
m_j\fp[G/G_j]$.  Now $M'$ may be taken to be $M'_n$.  
%\qed
\enddemo

%The data in Tables 2--6 of Section~8 were obtained using the
%algorithm described above.  
The data in Tables 8.2--8.6 of the next section were obtained using the
algorithm described above.

\head 8. {\it Tables and final remarks }\endhead
For each $p$, $n$, and $d$, let $\tilde K^*_{n,d}$ be the direct
summand of $K^*_{n,d}$ corresponding to the reduced Morava K-theory
$\tilde \K(B(\fp)^d)$.  Thus $\tilde K^k_{n,d}=K^k_{n,d}$ for $k\neq
0$, and $K^0_{n,d} = \tilde K^0_{n,d}\oplus T$, where $T$ is the trivial
$\fp[GL(V)]$-submodule of dimension one spanned by the monomial 1.
The $\fp$-dimension of $\tilde K^k_{n,d}$ is $(p^{nd}-1)/(p^n-1)$.  
Table~8.1 describes the $SL_2(\fthree)$-module structure of $\tilde
K^k_{n,2}$ (for $p=3$) in terms of the indecomposable modules
$I_1,\ldots,I_7$ as described in Section~6.

\midinsert
%Table 1.  
\centerline{\it Table 8.1: 
The $SL_2(\fthree)$-summands of $\tilde K^k_{n,2}$.} 
$$
\vbox{\tabskip=0pt \offinterlineskip
\def\tablerule{\noalign{\hrule}}

\halign{
  \bigstrut#&\strut#\tabskip=1em plus2em
  &\hfil#\hfil&\strut#&
  \hfil#\hfil&\strut#&
  \hfil#\hfil&\strut#&
  \hfil#\hfil&\strut#&
  \hfil#\hfil&\strut#&
  \hfil#\hfil&\strut#&
  \hfil#\hfil&\strut#&
  \hfil#\hfil&\strut#&
\hfil#\hfil&\strut#\tabskip=0pt\crcr
\tablerule
&\vrule&$n$&\vrule&$k$&\vrule&$I_1$&\vrule&$I_2$&
\vrule&$I_3$&\vrule&$I_4$&\vrule&$I_5$&\vrule&$I_6$&\vrule&$I_7$&\vrule\cr
\tablerule
&\vrule&1&\vrule&0&\vrule&1&\vrule&0&
\vrule&0&\vrule&0&\vrule&0&\vrule&0&\vrule&1&\vrule\cr
\tablerule
&\vrule&1&\vrule&1&\vrule&0&\vrule&0&
\vrule&0&\vrule&0&\vrule&1&\vrule&0&\vrule&0&\vrule\cr
\tablerule
&\vrule&2&\vrule&0,4&\vrule&1&\vrule&0&
\vrule&1&\vrule&0&\vrule&0&\vrule&0&\vrule&2&\vrule\cr
\tablerule
&\vrule&2&\vrule&1,3,5,7&\vrule&0&\vrule&0&
\vrule&0&\vrule&0&\vrule&1&\vrule&1&\vrule&0&\vrule\cr
\tablerule
&\vrule&2&\vrule&2,6&\vrule&1&\vrule&0&
\vrule&0&\vrule&0&\vrule&0&\vrule&0&\vrule&3&\vrule\cr
\tablerule
&\vrule&3&\vrule&even&\vrule&1&\vrule&0&
\vrule&2&\vrule&0&\vrule&0&\vrule&0&\vrule&7&\vrule\cr
\tablerule
&\vrule&3&\vrule&odd&\vrule&0&\vrule&0&
\vrule&0&\vrule&0&\vrule&1&\vrule&4&\vrule&0&\vrule\cr
\tablerule}}
$$
%\centerline{\it Table 1: 
%The $SL_2(\fthree)$-summands of $\tilde K^k_{n,2}$.} 
\endinsert

Let $V$ have dimension $d=3$ over $\fp$.  Let $l$ be a line in $V$,
and let $\pi$ be a plane in $V$ containing $l$.  The group $GL(V)$
acts on the set of all such pairs, and the stabilizer of the pair 
$(l,\pi)$ contains a unique Sylow $p$-subgroup $U(V)$ of $GL(V)$ 
(and is in fact equal to the normalizer of $U(V)$).  Let $C$ be a
generator for the centre of $U(V)$, which is cyclic of order $p$.  Let
$A$ be a non-central element of $U(V)$ stabilizing every line in
$\pi$, and let $B$ be a non-central element of $U(V)$ stabilizing
every plane containing $l$.  Then $A$ and $B$ generate $U(V)$, and
after replacing $C$ by a power if necessary, the commutator of $A$ and
$B$ is equal to $C$.  If we identify $V$ with $(\fp)^3$, and take $U(V)$ to be
the upper triangular matrices, then we may take
$$A= \pmatrix
1&1&0\\
0&1&0\\
0&0&1\endpmatrix
\qquad
B = \pmatrix
1&0&0\\
0&1&1\\
0&0&1\endpmatrix.$$
(Recall that $V$ is to be viewed as the space of row vectors with a
right $GL(V)$-action.)

\midinsert
\centerline{\it Table 8.2: 
The $D_8$-summands of $\tilde K^k_{n,3}$.} 
$$
\vbox{\tabskip=0pt \offinterlineskip
\def\tablerule{\noalign{\hrule}}

\halign{
  \bigstrut#&\strut#\tabskip=1em plus2em
  &\hfil#\hfil&\strut#&
  \hfil#\hfil&\strut#&
  \hfil#\hfil&\strut#&
  \hfil#\hfil&\strut#&
  \hfil#\hfil&\strut#&
  \hfil#\hfil&\strut#&
  \hfil#\hfil&\strut#&
  \hfil#\hfil&\strut#&
  \hfil#\hfil&\strut#&
\hfil#\hfil&\strut#\tabskip=0pt\crcr
\tablerule
&\vrule&$n$&\vrule&$k$&\vrule&$P_1$&\vrule&$P_2$&
\vrule&$P_3$&\vrule&$P_4$&\vrule&$P_5$&\vrule&$P_6$
&\vrule&$P_7$&\vrule&$P_8$&\vrule\cr
\tablerule
&\vrule&1&\vrule&0&\vrule&0&\vrule&1&\vrule&0&
\vrule&0&\vrule&0&\vrule&0&\vrule&1&\vrule&1&\vrule\cr
\tablerule
&\vrule&2&\vrule&0&\vrule&0&\vrule&2&\vrule&2&
\vrule&1&\vrule&0&\vrule&0&\vrule&0&\vrule&1&\vrule\cr
\tablerule
&\vrule&2&\vrule&1,2&\vrule&1&\vrule&1&\vrule&1&
\vrule&0&\vrule&0&\vrule&1&\vrule&1&\vrule&1&\vrule\cr
\tablerule
&\vrule&3&\vrule&0&\vrule&4&\vrule&4&\vrule&4&
\vrule&2&\vrule&0&\vrule&0&\vrule&0&\vrule&1&\vrule\cr
\tablerule
&\vrule&3&\vrule&1,6&\vrule&6&\vrule&1&\vrule&3&
\vrule&0&\vrule&0&\vrule&3&\vrule&1&\vrule&1&\vrule\cr
\tablerule
&\vrule&3&\vrule&2,5&\vrule&5&\vrule&3&\vrule&3&
\vrule&1&\vrule&0&\vrule&1&\vrule&1&\vrule&1&\vrule\cr
\tablerule
&\vrule&3&\vrule&3,4&\vrule&5&\vrule&2&\vrule&4&
\vrule&1&\vrule&0&\vrule&2&\vrule&0&\vrule&1&\vrule\cr
\tablerule
&\vrule&4&\vrule&0&\vrule&24&\vrule&8&\vrule&8&
\vrule&4&\vrule&0&\vrule&0&\vrule&0&\vrule&1&\vrule\cr
\tablerule
&\vrule&4&\vrule&1,14&\vrule&28&\vrule&1&\vrule&7&
\vrule&0&\vrule&0&\vrule&7&\vrule&1&\vrule&1&\vrule\cr
\tablerule
&\vrule&4&\vrule&2,13&\vrule&25&\vrule&7&\vrule&7&
\vrule&3&\vrule&0&\vrule&1&\vrule&1&\vrule&1&\vrule\cr
\tablerule
&\vrule&4&\vrule&3,12&\vrule&27&\vrule&2&\vrule&8&
\vrule&1&\vrule&0&\vrule&6&\vrule&0&\vrule&1&\vrule\cr
\tablerule
&\vrule&4&\vrule&4,11&\vrule&25&\vrule&6&\vrule&8&
\vrule&3&\vrule&0&\vrule&2&\vrule&0&\vrule&1&\vrule\cr
\tablerule
&\vrule&4&\vrule&5,10&\vrule&27&\vrule&3&\vrule&7&
\vrule&1&\vrule&0&\vrule&5&\vrule&1&\vrule&1&\vrule\cr
\tablerule
&\vrule&4&\vrule&6,9&\vrule&26&\vrule&5&\vrule&7&
\vrule&2&\vrule&0&\vrule&3&\vrule&1&\vrule&1&\vrule\cr
\tablerule
&\vrule&4&\vrule&7,8&\vrule&26&\vrule&4&\vrule&8&
\vrule&2&\vrule&0&\vrule&4&\vrule&0&\vrule&1&\vrule\cr
\tablerule
&\vrule&5&\vrule&0&\vrule&112&\vrule&16&\vrule&16&
\vrule&8&\vrule&0&\vrule&0&\vrule&0&\vrule&1&\vrule\cr
\tablerule}}
$$
%\centerline{\it Table 2: 
%The $D_8$-summands of $\tilde K^k_{n,3}$.} 
\endinsert

\midinsert
\centerline{\it Table 8.3: 
A maximal $\hbox{\rm Syl}_3(GL_3(\fthree))$-permutation submodule 
of $\tilde K^k_{n,3}$.} 
$$
\vbox{\tabskip=0pt \offinterlineskip
\def\tablerule{\noalign{\hrule}}

\halign{
  \bigstrut#&\strut#\tabskip=1em plus2em
  &\hfil#\hfil&\strut#&
  \hfil#\hfil&\strut#&
  \hfil#\hfil&\strut#&
  \hfil#\hfil&\strut#&
  \hfil#\hfil&\strut#&
  \hfil#\hfil&\strut#&
  \hfil#\hfil&\strut#&
  \hfil#\hfil&\strut#&
  \hfil#\hfil&\strut#&
\hfil#\hfil&\strut#\tabskip=0pt\crcr
\tablerule
&\vrule&$n$&\vrule&$k$&\vrule&dim.~$M'$&\vrule&$P_1$&
\vrule&$P_2$&\vrule&$P_5$&\vrule&$P_6$&\vrule&$P_7$
&\vrule&$P_{10}$&\vrule&$P_{11}$&\vrule\cr
\tablerule
&\vrule&1&\vrule&0,1&\vrule&13&\vrule&0&\vrule&1&
\vrule&0&\vrule&0&\vrule&0&\vrule&1&\vrule&1&\vrule\cr
\tablerule
&\vrule&2&\vrule&0&\vrule&91&\vrule&1&\vrule&3&
\vrule&3&\vrule&1&\vrule&0&\vrule&0&\vrule&1&\vrule\cr
\tablerule
&\vrule&2&\vrule&1&\vrule&65&\vrule&1&\vrule&1&
\vrule&2&\vrule&0&\vrule&3&\vrule&0&\vrule&2&\vrule\cr
\tablerule
&\vrule&2&\vrule&2&\vrule&71&\vrule&1&\vrule&2&
\vrule&2&\vrule&0&\vrule&1&\vrule&1&\vrule&2&\vrule\cr
\tablerule
&\vrule&2&\vrule&3&\vrule&73&\vrule&1&\vrule&2&
\vrule&2&\vrule&0&\vrule&1&\vrule&2&\vrule&1&\vrule\cr
\tablerule
&\vrule&2&\vrule&4&\vrule&57&\vrule&1&\vrule&1&
\vrule&1&\vrule&0&\vrule&2&\vrule&1&\vrule&3&\vrule\cr
\tablerule
&\vrule&2&\vrule&5&\vrule&65&\vrule&1&\vrule&2&
\vrule&1&\vrule&0&\vrule&1&\vrule&2&\vrule&2&\vrule\cr
\tablerule
&\vrule&2&\vrule&6&\vrule&67&\vrule&1&\vrule&2&
\vrule&1&\vrule&0&\vrule&2&\vrule&2&\vrule&1&\vrule\cr
\tablerule
&\vrule&2&\vrule&7&\vrule&73&\vrule&1&\vrule&2&
\vrule&2&\vrule&0&\vrule&2&\vrule&1&\vrule&1&\vrule\cr
\tablerule}}
$$
%\centerline{\it Table 3: 
%A maximal $\hbox{\rm Syl}_3(GL_3(\fthree))$-permutation submodule 
%of $\tilde K^k_{n,3}$.} 
\endinsert

For $p=2$ the group $U(V)$ has 8 conjugacy classes of subgroups, which
we list in the following order:  
$$\{1\},\langle A\rangle, \langle B\rangle, \langle C\rangle, 
\langle AB\rangle, \langle A,C\rangle, \langle B,C\rangle, U(V). 
$$
Let $P_1,\ldots,P_8$ be the corresponding transitive permutation 
modules, so that $P_1$ is the free module and $P_8$ is the trivial
module.  Similarly, for $p>2$, $U(V)$ has $2p+5$ conjugacy classes of
subgroups, which we list as:  
$$\eqalign{\{1\}, 
\langle A \rangle, \langle AB \rangle&,\ldots, \langle AB^{p-1} \rangle,  
\langle B \rangle, \langle C \rangle,\cr &\langle A,C \rangle, 
\langle AB,C \rangle,\ldots, \langle AB^{p-1},C \rangle, 
\langle B,C \rangle, U(V).}
$$
Again we let $P_1,\ldots, P_{2p+5}$ be the corresponding transitive
permutation modules.

\midinsert
\centerline{\it Table 8.4: 
A maximal $\hbox{\rm Syl}_5(GL_3(\ffive))$-permutation submodule 
of $\tilde K^k_{n,3}$.} 
$$
\vbox{\tabskip=0pt \offinterlineskip
\def\tablerule{\noalign{\hrule}}

\halign{
  \bigstrut#&\strut#\tabskip=1em plus2em
  &\hfil#\hfil&\strut#&
  \hfil#\hfil&\strut#&
  \hfil#\hfil&\strut#&
  \hfil#\hfil&\strut#&
  \hfil#\hfil&\strut#&
  \hfil#\hfil&\strut#&
  \hfil#\hfil&\strut#&
  \hfil#\hfil&\strut#&
  \hfil#\hfil&\strut#&
\hfil#\hfil&\strut#\tabskip=0pt\crcr
\tablerule
&\vrule&$n$&\vrule&$k$&\vrule&dim.~$M'$&\vrule&$P_1$&
\vrule&$P_2$&\vrule&$P_7$&\vrule&$P_8$&\vrule&$P_9$
&\vrule&$P_{14}$&\vrule&$P_{15}$&\vrule\cr
\tablerule
&\vrule&1&\vrule&0--3&\vrule&31&\vrule&0&\vrule&1&
\vrule&0&\vrule&0&\vrule&0&\vrule&1&\vrule&1&\vrule\cr
\tablerule
&\vrule&2&\vrule&0&\vrule&651&\vrule&3&\vrule&5&
\vrule&5&\vrule&1&\vrule&0&\vrule&0&\vrule&1&\vrule\cr
\tablerule
&\vrule&2&\vrule&1&\vrule&527&\vrule&3&\vrule&1&
\vrule&4&\vrule&0&\vrule&5&\vrule&0&\vrule&2&\vrule\cr
\tablerule
&\vrule&2&\vrule&2&\vrule&447&\vrule&2&\vrule&3&
\vrule&4&\vrule&0&\vrule&3&\vrule&1&\vrule&2&\vrule\cr
\tablerule
&\vrule&2&\vrule&3&\vrule&467&\vrule&2&\vrule&4&
\vrule&4&\vrule&0&\vrule&2&\vrule&1&\vrule&2&\vrule\cr
\tablerule
&\vrule&2&\vrule&4&\vrule&587&\vrule&3&\vrule&4&
\vrule&4&\vrule&0&\vrule&1&\vrule&1&\vrule&2&\vrule\cr
\tablerule
&\vrule&2&\vrule&5&\vrule&591&\vrule&3&\vrule&4&
\vrule&4&\vrule&0&\vrule&1&\vrule&2&\vrule&1&\vrule\cr
\tablerule}}
$$
%\centerline{\it Table 4: 
%A maximal $\hbox{\rm Syl}_5(GL_3(\ffive))$-permutation submodule 
%of $\tilde K^k_{n,3}$.} 
\endinsert

Tables 8.2, 8.3, and 8.4 describe maximal $U(V)$-permutation submodules $M'$
of $\tilde K^k_{n,3}$ in the cases $p=2$, 3, and 5 respectively.
These submodules were found using the algorithm of Proposition~7.1,
with the conjugacy classes of subgroups of $U(V)$ listed in the order
given above.  The permutation modules omitted from Tables 8.3~and~8.4
never arose as summands of any such $M'$.  In Table~8.2 the dimension of
$M'$ is omitted since in these cases $M'$ was always the whole of
$\tilde K^k_{n,3}$, with dimension $2^{2n}+2^n+1$.  The dimensions of 
$\tilde K^k_{1,3}$ and $\tilde K^k_{2,3}$ are 13 and 91 respectively for
$p=3$, and 31 and 651 respectively for $p=5$.  Table~8.4 in the
case $n=2$ is incomplete in the sense that not all values of $k$ have
been considered.  This is because each row required over 24 hours'
computing time.  

Note that for fixed $p$ and $n$, the modules $K^k_{n,3}$ tend
not to be isomorphic to each other, except when $n=1$, and isomorphic
in pairs when $p=2$.  For each $p$, the ring structure gives rise to a 
duality 
$$ K^k_{n,3}\times  K^{N-k}_{n,3}\rightarrow \fp \subseteq K^N_{n,3},$$
where $N=p^n-1$.  This explains the observed fact that whenever
$K^k_{n,3}$ is a permutation module, then 
$K^{N-k}_{n,3}\cong K^k_{n,3}$.  

We find it intriguing that in the case $p=2$ we have been unable to
find pairs $(n,k)$ such that $K^k_{n,3}$ is not a $U(V)$-permutation 
module.  Note also that for each $n$ and $p$ considered, $K^0_{n,3}$
is a $U(V)$-permutation module, although it is easy to show that
usually $K^0_{n,3}$ cannot be a $GL(V)$-permutation module by
comparing the information in the tables with the information given by
Brauer characters.  (This technique may be used to prove the cases 
$n=2$, 4~and~5 of Theorem~1.1(d), which we leave as an exercise.)  

Finally, in Tables 8.5~and~8.6 we give just enough information to prove
Theorem~1.4, in the cases $p=3$ and $p=5$ respectively.  That is, for
each subgroup $H$ of $U(V)$ of order $p^2$, we give the dimension of a
maximal $H$-permutation submodule~$M''$ of~$K^1_{2,3}$.  The dimension
of $K^1_{2,3}$ is 91 for $p=3$ and 651 for $p=5$.  Only one of the
subgroups $\langle AB,C \rangle,\ldots, \langle AB^{p-1},C \rangle$ is
listed in these tables, because these subgroups are all conjugate in
$GL(V)$ and so give rise to $M''$'s of the same dimension.  The
programs were run separately for each of these groups however, as a
check.

\midinsert
\centerline{\it Table 8.5: A maximal $H$-permutation submodule of 
$K^1_{2,3}$ ($p=3$).}
$$
\vbox{\tabskip=0pt \offinterlineskip
\def\tablerule{\noalign{\hrule}}

\halign{
  \bigstrut#&\strut#\tabskip=1em plus2em
  &\hfil#\hfil&\strut#&
\hfil#\hfil&\strut#\tabskip=0pt\crcr
\tablerule
&\vrule&Subgroup $H$&\vrule&dim.~$M''$&\vrule\cr
\tablerule
&\vrule&$\langle A,C\rangle$&\vrule&69&\vrule\cr
\tablerule
&\vrule&$\langle AB,C\rangle$&\vrule&84&\vrule\cr
\tablerule
&\vrule&$\langle B,C\rangle$&\vrule&87&\vrule\cr
\tablerule}}
$$
%\centerline{\it Table 5: A maximal $H$-permutation submodule of 
%$K^1_{2,3}$ ($p=3$).}

\endinsert

\midinsert
\centerline{\it Table 8.6: A maximal $H$-permutation submodule of
$K^1_{2,3}$ ($p=5$).}
$$
\vbox{\tabskip=0pt \offinterlineskip
\def\tablerule{\noalign{\hrule}}

\halign{
  \bigstrut#&\strut#\tabskip=1em plus2em
  &\hfil#\hfil&\strut#&
\hfil#\hfil&\strut#\tabskip=0pt\crcr
\tablerule
&\vrule&Subgroup $H$&\vrule&dim.~$M''$&\vrule\cr
\tablerule
&\vrule&$\langle A,C\rangle$&\vrule&535&\vrule\cr
\tablerule
&\vrule&$\langle AB,C\rangle$&\vrule&628&\vrule\cr
\tablerule
&\vrule&$\langle B,C\rangle$&\vrule&643&\vrule\cr
\tablerule}}
$$
%\centerline{\it Table 6: A maximal $H$-permutation submodule of
%$K^1_{2,3}$ ($p=5$).}

\endinsert

Kriz's example of a 3-group $G$ such that $K(2)^*(BG)$ is not
concentrated in even degrees is the Sylow 3-subgroup of
$GL_4(\fthree)$ [\kri].  This group is expressible as the split extension
with kernel $(\fthree)^3$ and quotient the Sylow 3-subgroup of
$GL_3(\fthree)$, with the natural action.  There may be a connection
between our result that $K^*_{2,3}$ is not a permutation module for 
the Sylow 3-subgroup of $GL_3(\fthree)$ and the fact that $K(2)^*(BG)$
is not entirely even.  If so, then Theorem~1.4 suggests if $H$ is a split
extension with kernel $(\fthree)^3$ and quotient a subgroup of
$GL_3(\fthree)$ of order nine, then possibly $K(2)^*(BH)$ is not
entirely even.  Such $H$ include the extraspecial group of order $3^5$
and exponent~3.  
\medskip
\par\noindent
{\it Acknowledgements. }Throughout this project the
second-named author was supported by an E.C. Leibniz Fellowship at the
CRM{}.  The first-named author was initially supported by a DGICYT
Fellowship at the CRM and later by an E.C. Leibniz Fellowship at the
MPI{}.  Our computer programs were run on a network of Sun SPARC
computers at the MPI{}.  The authors gratefully acknowledge the
hospitality of both the CRM and the MPI{}.

%%%%%%%%%%%%%%%%%%%%%%%%%%%%%%%%%%%%%%%%%%%%%%%%%%%%%%%%%%%%%%%%%%%%%%%
%                                                                     %
%       References for the K(n)^*(BV)-paper                           %
%								      %
%%%%%%%%%%%%%%%%%%%%%%%%%%%%%%%%%%%%%%%%%%%%%%%%%%%%%%%%%%%%%%%%%%%%%%%

\bigskip
\centerline{\smc References}
{\sm
\advance\baselineskip by -2pt
\roster
\item"{\sm [{\smb 1}]}"
{\smallcaps J.~L.~Alperin}. {\smit Local representation theory}
(Cambridge University Press, 1986).

\item"{\sm [{\smb 2}]}"
{\smallcaps M. Brunetti}. On the canonical $\scriptstyle GL_2(\ftwo)$-module 
structure of $\scriptstyle K(n)^*(B\zz/2\times B\zz/2)$. (To appear.)

\item"{\sm [{\smb 3}]}"
{\smallcaps C.~W.~Curtis} and {\smallcaps I.~Reiner}
{\smit Methods of Representation theory I} 
(Wiley, 1981).

\item"{\sm [{\smb 4}]}"
{\smallcaps M.~J.~Hopkins, N.~J.~Kuhn} and {\smallcaps D.~C.~Ravenel}.
Generalized Group Characters and Complex Oriented Cohomology Theories.
Preprint (1992).

\item"{\sm [{\smb 5}]}"
{\smallcaps J.~R.~Hunton}. The Morava K-theory of wreath products.
{\smit Math. Proc. Camb. Phil. Soc.} {\smb 107} (1990), 309-318.

\item"{\sm [{\smb 6}]}"
{\smallcaps I.~Kriz}. Morava K-theories of classifying spaces: some 
calculations. Preprint (1995).

\item"{\sm [{\smb 7}]}"
{\smallcaps N.~J.~Kuhn}. Morava K-theories of some classifying spaces. 
{\smit Trans. Amer. Math. Soc.} {\smb 304} (1987), 193-205.

\item"{\sm [{\smb 8}]}"
{\smallcaps N.~J.~Kuhn}. The mod {\smit p} K-theory of classifying spaces of finite 
groups. {\smit Jour. Pure and Applied Algebra} {\smb 44} (1987), 269-271.

\item"{\sm [{\smb 9}]}"
{\smallcaps D.~C.~Ravenel} and {\smallcaps W.~S.~Wilson}. The Morava K-theories 
of Eilenberg-Mac Lane spaces and the Conner-Floyd conjecture.
{\smit Amer. Jour. Math.} {\smb 102} (1980),  691-748.

\item"{\sm [{\smb 10}]}"
{\smallcaps M.~Sch\"onert} et al.
{\smit GAP (Groups, Algorithms and Programming) Version 3  Release 2}
(RWTH Aachen, 1993).

\item"{\sm [{\smb 11}]}"
{\smallcaps B.~Schuster}. On the Morava K-theory of some finite 2-groups.
Preprint (1994).

\item"{\sm [{\smb 12}]}"
{\smallcaps M. Tanabe}. On Morava K-theories of Chevalley groups.
{\smit Amer. Jour. of Math.} {\smb 117} (1995), 263-278.

\item"{\sm [{\smb 13}]}" % Tezuka-Yagita 1, tex doesn't like keys like \teya1
{\smallcaps M. Tezuka, N. Yagita}. Cohomology of finite groups and 
Brown-Peterson cohomology. In {\smit Algebraic Topology} (Arcata, CA, 1986)
(Springer, 1989), 396-408.

\item"{\sm [{\smb 14}]}"
{\smallcaps M. Tezuka, N. Yagita}. Cohomology of groups and Brown-Peterson 
cohomology II. In {\smit Homotopy theory and related topics (Kinosaki, 1988)}
(Springer, 1990), 57-69.

\item"{\sm [{\smb 15}]}"
{\smallcaps W.~S.~Wilson}. Brown-Peterson Homology, an Introduction and 
Sampler.
{\smit Regional Conference Series in Mathematics, No. {\smb 48}}, (AMS, 1980).

\item"{\sm [{\smb 16}]}"
{\smallcaps C.~Wilkerson}. A primer on Dickson Invariants.
{\smit Proceedings of the Northwestern Homotopy conference},
A.M.S. Contemp. Math. Series {\smb 19} (1983), 421-434.

\item"{\sm [{\smb 17}]}"
{\smallcaps N. Yagita}. Note on BP-theory for extensions of cyclic groups by 
elementary abelian p-groups.
Preprint (1994).

\endroster
\bigskip
\bigskip
\bigskip
%\enddocument
%\end

$$
\hbox{\vbox{
\hbox{Ian Leary}
\hbox{Max-Planck Institut f\"ur Mathematik}
%\hbox{Gottfried Claren Str.\ 26}
\hbox{53225 Bonn}
\hbox{Germany}
%\hbox{email: leary\@mpim-bonn.mpg.de}
\hbox{\sm From Jan.\ 1996:}
\hbox{Faculty of Math. Studies}
\hbox{Univ. of Southampton}
\hbox{Southampton SO17 1BJ}
\hbox{England}
}\qquad
\vbox{\hbox{Bj\"orn Schuster}
\hbox{Centre de Recerca Matem\`atica}
%\hbox{Institut d'Estudis Catalans}
\hbox{E--08193 Bellaterra (Barcelona)}
\hbox{Spain}
%\hbox{email: schuster\@bianya.crm.es}
%\hbox{\sm From Sept.\ 1995:} 
%\hbox{The Fields Institute}
%\hbox{222 College Street}
%\hbox{Toronto, Ontario}
%\hbox{Canada, M5T 3J1}}}
\hbox{\sm From June 1996:}
\hbox{FB Mathematik}
\hbox{Bergische Universit\"at}
\hbox{D--42097 Wuppertal}
\hbox{Germany}}}
$$

\enddocument